\documentclass[12pt]{article}
\usepackage{amsfonts,amssymb,amsbsy,amsmath,amsthm,times}
\usepackage{endnotes}
\newtheorem{theorem}{Theorem}[section]

\newtheorem{lemma}[theorem]{Lemma}
\newtheorem{follow}[theorem]{Corollary}

\newtheorem{pr}[theorem]{Proposition}
\theoremstyle{definition}

\newcommand{\bel}{\begin{equation} \label}
\newcommand{\ee}{\end{equation}}

\newcommand{\xp}{x}
\newcommand{\mper}{m_\perp}

\newcommand{\rt}{{\mathbb R}^{3}}
\newcommand{\rd}{{\mathbb R}^{2}}
\newcommand{\re}{{\mathbb R}}

\newcommand{\cd}{{\mathbb C}^{2}}

\newcommand{\phe}{\tilde{\varphi}}

{

\begin{document}

\begin{center}

{\large \bf Low Energy Asymptotics of the Spectral Shift Function for Pauli Operators
with Nonconstant Magnetic Fields}

\let\thefootnote\relax\footnotetext{AMS 2010 Mathematics Subject Classification: 35J10, 81Q10, 35P20,
35P25, 47F05}

\let\thefootnote\relax\footnotetext{Keywords:
Pauli operators, spectral shift function, low-energy asymptotics,
Levinson formula}

\end{center}

\vspace{0.5cm}

\begin{center}
{\sc Georgi  D. Raikov}*
\let\thefootnote\relax\footnotetext{*Facultad de Matem\'aticas,
Pontificia Universidad Cat\'olica de Chile, Vicu\~na Mackenna 4860,}
\let\thefootnote\relax\footnotetext{Macul, Santiago de Chile.
E-mail: graikov@mat.puc.cl}  \\

\end{center}

\medskip

\begin{abstract}
We consider the 3D Pauli operator with nonconstant magnetic field
${\bf B}$ of constant direction, perturbed by a symmetric
matrix-valued electric potential $V$ whose coefficients decay fast
enough at infinity. We investigate the low-energy asymptotics of
the corresponding spectral shift function. As a corollary, for
generic negative $V$, we obtain a generalized Levinson formula,
relating the low-energy asymptotics of the eigenvalue
counting function and of the scattering phase of the perturbed operator.\\

\end{abstract}

\section{Introduction} \label{s1} \setcounter{equation}{0}
Suppose that the magnetic field ${\bf B} : \rt \to \rt$ has a
constant direction, say,
    \bel{sn50}
{\bf B} = (0,0,b).
    \ee
By the Maxwell equation, ${\rm div}\,{\bf B} = 0$, we should then
have $\frac{\partial b}{\partial x_3} = 0$. Assume that the
function $b : \rd \to \re$ is continuous and bounded. In
Subsection \ref{ss20} we describe in more detail the class of
admissible functions $b$. Let ${\bf A} \in C^1(\rt; \rt)$ be a
magnetic potential generating the magnetic field ${\bf B}$, i.e.
${\rm curl}\,{\bf A} = {\bf B}$. Introduce the Pauli matrices
$$
\hat{\sigma}_1 = \left(
\begin{array} {cc}
0 & 1\\
1 & 0
\end{array}
\right), \quad
\hat{\sigma}_2 = \left(
\begin{array} {cc}
0 & -i\\
i & 0
\end{array}
\right), \quad
\hat{\sigma}_3  = \left(
\begin{array} {cc}
1 & 0\\
0 & -1
\end{array}
\right).
$$
 Set $\hat{\sigma} : = (\hat{\sigma}_1,
\hat{\sigma}_2, \hat{\sigma}_3)$. Let
    \bel{sn61}
H_0 : = \left(\hat{\sigma} \cdot (-i\nabla - {\bf A})\right)^2
    \ee
be the unperturbed self-adjoint Pauli operator defined originally on $C_0^{\infty}(\re^3; {\mathbb C}^2)$,
and then closed in $L^2(\re^3; {\mathbb C}^2)$. We have
$$
H_0 : = \left(
\begin{array} {cc}
(-i\nabla - {\bf A})^2 - b & 0\\
0 & (-i\nabla - {\bf A})^2 + b
\end{array}
\right) : =
\left(
\begin{array} {cc}
H_0^- & 0\\
0 & H_0^+
\end{array}
\right) = H_0^- \oplus H_0^+.
$$
Further, let $v_{jk} \in L^{\infty}(\re^3)$, $j,k = 1,2$. Assume
that  $v_{11}$ and $v_{22}$ are real-valued, and $v_{12} =
\overline{v_{21}}$. Introduce the symmetric  matrix
$$
V({\bf x}) : = \left(
\begin{array} {cc}
v_{11}({\bf x}) & v_{12}({\bf x})\\
v_{21}({\bf x}) & v_{22}({\bf x})
\end{array}
\right), \quad {\bf x} \in \re^3.
$$
On the domain of $H_0$ define the operator
 $$
H : = H_0 + V.
$$
Assume that
    \bel{13}
    (H - i)^{-1} - (H_0 - i)^{-1} \in S_1(L^2(\re^3; {\mathbb
    C}^2))
    \ee
    where $S_1(X)$ denotes the trace class of linear operators acting in the Hilbert space $X$ .
    By the diamagnetic inequality and the
    boundedness of $b$, we find that \eqref{13} holds true if
    \bel{14}
    |v_{jk}|^{1/2} (-\Delta + 1)^{-1} \in S_2(L^2(\re^3)), \quad j,k = 1,2,
    \ee
    where $S_2(X)$ denotes the Hilbert-Schmidt class of linear operators acting in  $X$.
    On its turn, \eqref{14} holds true if and only if $v_{jk} \in L^{1}(\re^3)$.\\
By \eqref{13}, there exists a unique  $\xi = \xi(\cdot;H,H_0)
 \in L^1(\re; (1+E^2)^{-1}dE)$ which vanishes identically on
 $(-\infty,\inf \sigma(H))$, such that {\em the Lifshits-Krein trace formula}
    \bel{sn15}
{\rm Tr}\,(f(H)-f(H_0)) = \int_{\re} \xi(E; H,H_0) f'(E) dE
    \ee
    holds for each $f \in C^{\infty}_0(\re)$
(see the original works \cite{lif}, \cite{kr},  or \cite[Chapter
8]{y}).\\
 The function $\xi(\cdot;H,H_0)$ is called
{\em the spectral shift function} (SSF) for the operator pair
$(H,H_0)$. If $E < 0 = \inf \sigma(H_0)$, then the spectrum of $H$
below $E$ could be at most discrete, and for almost every $E<0$ we
have
\begin{equation} \label{101a}
\xi(E; H, H_0) = - N(E; H)
\end{equation}

 where   $ N(E; H)$ denotes the number of
eigenvalues of $H$ lying in the interval $(-\infty,E)$, and
counted with their multiplicities.   On the other hand, for almost
every $E \in \sigma_{\rm ac}(H_0) = [0,\infty)$ (see Corollary
\ref{rtf21a} below), the SSF $\xi(E; H, H_0)$ is related to the
scattering determinant ${\rm det}\;S(E; H, H_0)$ for the pair $(H,
H_0)$ by {\em the Birman-Krein formula}
    \bel{sn30}
{\rm det}\;S(E; H, H_0) = e^{-2\pi i \xi(E; H, H_0)}
    \ee
(see the original work \cite{bk} or \cite[Section 8.4]{y}). \\
A priori, the SSF $\xi(E;H,H_0)$ is defined for almost every $E
\in \re$. In this article, if $E \in (-\infty, {\mathcal
C})\setminus\{0\}$ where ${\mathcal C}>0$ is a constant defined in
\eqref{22}, we will identify $\xi(E;H,H_0)$  with a representative
of its equivalence class, described explicitly below in Subsection
\ref{ss31} under the assumption that the matrix $V({\bf x})$,
${\bf x} \in \re^3$, has a definite sign. Under our generic
assumptions on $V$, we check that the SSF $\xi(\cdot;H,H_0)$ is
bounded on every compact subset of $(-\infty,
\mathcal{C})\setminus\{0\}$, and continuous on $(-\infty,
\mathcal{C})\setminus(\{0\} \cup \sigma_{\rm pp}(H))$ where
$\sigma_{\rm pp}(H)$ denotes the set of the
eigenvalues of $H$ (see Proposition \ref{rtf21} below). \\
The main results of the article concern the asymptotic behavior of
the SSF $\xi(E;H,H_0)$ as $E \to 0$ for perturbations $V$ of
definite sign. We show that even for certain $V$ of  compact
support, the SSF $\xi(\cdot;H,H_0)$ has a singularity at the
origin (see Theorems \ref{frt1} and \ref{frt2} below). More
precisely, we show that $\xi(E;H,H_0) \to \infty$ as $E \downarrow
0$ if the perturbation is positive, and $\xi(E;H,H_0) \to -\infty$
as $E \uparrow 0$ and $E \downarrow 0$ if the perturbation is
negative. The singularities of the SSF at the origin are described
in the terms of effective Hamiltonians of Berezin-Toeplitz type;
their spectral properties have been studied, for instance, in
\cite{r0}, \cite{rw}, and \cite{r5}. Assuming that the
perturbation admits a power-like or exponential decay at infinity,
or that it has a compact support, we obtain the first asymptotic
term of $\xi(E;H,H_0)$ as $E \uparrow 0$ and $E \downarrow 0$ (see
Corollaries \ref{frf01} and \ref{frf02} below). In particular,  if
the perturbation is negative, we show that there exists a finite
positive limit
$$
\lim_{E \downarrow 0} \frac{\xi(E;H,H_0)}{\xi(-E;H,H_0)}
$$
which depends only on the decay rate of $V$ at infinity (see Corollary  \ref{frf03} below).\\
Similar results concerning the singularities at the Landau levels
of the SSF in case where the  unperturbed operator is the 3D
Schr\"odinger operator with constant magnetic field, and the
perturbation is a sign-definite scalar potential which decays fast
enough at infinity, were obtained in \cite{fr}. The relation
between these singularities and
the possible accumulation of resonances at the Landau levels, was considered in \cite{bbr}. \\
The paper is organized as follows. In Section \ref{s2} we discuss
the class of the admissible magnetic fields, describe the basic
spectral properties of the operator $H_0$, and introduce the
Berezin-Toeplitz operators we need. In Section \ref{s3} we
formulate our main results as well as some corollaries of them.
Section \ref{s4} is devoted to
 auxiliary material such as   the representation of the
SSF due to A. Pushnitski, and   estimates of appropriate
sandwiched resolvents.  Finally, Section \ref{s5} contains the
proofs of Theorems \ref{frt1} -- \ref{frt2}.

\section{Admissible Magnetic Fields and Effective Hamiltonians}
 \label{s2}
 \setcounter{equation}{0}
\subsection{Admissible magnetic fields}
    \label{ss20}
Let  ${\bf B}$ have the form \eqref{sn50}.
 Assume that  $b = b_0
+ \tilde{b}$ where $b_0 > 0$ is a constant, while the function
$\tilde{b} : \rd \to \re$ is such that the Poisson equation
     \bel{rt00}
\Delta \tilde{\varphi} = \tilde{b}
    \ee
admits a  solution
     $\tilde{\varphi} : \rd \to \re$, continuous and bounded together with its derivatives of order up
to two. Abusing slightly the terminology, we will say that $b$ is
{\em an admissible} magnetic field. Also, we will call the
constant $b_0$ {\em the mean value} of $b$, and $\tilde{b}$ {\em
the background} of $b$. In our leading example, the admissible
background $\tilde{b}$ has the form
    \bel{sn0}
\tilde{b}(x) = \int_{\rd} e^{i\lambda \cdot x} d\nu(\lambda),
\quad x \in \rd,
    \ee
 where $\nu$ is a Borel charge (i.e. a complex-valued
measure) defined on $\rd$ which satisfies
    \bel{sn33}
    |\nu|(\rd) < \infty,
    \ee
    \bel{sn34} \nu(\delta) = \overline{\nu(-\delta)}
    \ee
     for each Borel set $\delta \subset \rd$,
     \bel{sn35}
    \nu(\left\{0\right\}) = 0,
    \ee
    and
    \bel{sn36}
    \int_{\rd} |\lambda|^{-2} d|\nu|(\lambda) < \infty.
    \ee
If $\tilde{b}$ satisfies \eqref{sn0}, then the Poisson equation
    \eqref{rt00} admits a solution
    \bel{sn1}
 \tilde{\varphi}(x)  : = - \int_{\rd} |\lambda|^{-2} e^{i\lambda \cdot x} d\nu(\lambda), \quad x \in \rd,
    \ee
which possesses all the prescribed properties.\\
Let us give two further examples of admissible backgrounds $\tilde{b}$ of the form \eqref{sn0}.\\
(i) Let $\lambda_n \subset \rd\setminus\{0\}$, $b_n \in {\mathbb
C}$, $n \in {\mathbb N}$. Assume that $\sum_{n \in {\mathbb N}}
|b_n| (1 + |\lambda_n|^{-2})< \infty$. Then the almost periodic
function $\tilde{b}(x) : = \sum_{n \in {\mathbb N}} b_n
e^{i\lambda_n . x}$, $x \in \rd$, is an admissible background,
provided that it is real-valued. In this case the charge $\nu$ in
\eqref{sn0} is singular with respect to the Lesbegue measure in
$\rd$. Evidently, the real-valued periodic functions with zero
mean value and absolutely convergent series of the Fourier
coefficients, belong to the described class of
admissible backgrounds.  \\
(ii) Let $f: \rd \to {\mathbb C}$ be a Lebesgue measurable
function which satisfies $f(\lambda) = \overline{f(-\lambda)}$,
$\lambda \in \rd$, and $\int_{\rd} (1 + |\lambda|^{-2})
|f(\lambda)| d\lambda < \infty$. Then $\tilde{b}(x) : = \int_{\rd}
e^{i\lambda \cdot x} f(\lambda) d\lambda$ is again an admissible
background. In this case charge $\nu$ in \eqref{sn0} is absolutely
continuous with respect to the Lesbegue
measure in $\rd$.\\
For $(x_1,x_2) \in \rd$ set $\varphi_0 : = b_0(x_1^2 + x_2^2)/4$
and
    \bel{sn2}
\varphi :  = \varphi_0 + \tilde{\varphi},
    \ee
$\tilde{\varphi}$ being introduced in \eqref{rt00}. Then $\Delta
\varphi_0 = b_0$ and $\Delta \varphi = b$. Put ${\bf A}: = (A_1,
A_2, A_3)$ with
    \bel{sn62}
A_1 : = - \frac{\partial \varphi}{\partial  x_2}, \quad A_2 : =
\frac{\partial \varphi}{\partial  x_1}, \quad A_3 = 0.
    \ee
The magnetic potential ${\bf A} : = (A_1,A_2,A_3) \in
C^1(\rt,\rt)$ generates the magnetic field ${\bf B}  = {\rm
curl}\,{\bf A} = (0,0,b)$. Changing, if necessary, the gauge, we
will assume that the magnetic potential ${\bf A}$ in \eqref{sn61}
is given by \eqref{sn62}.
    \subsection{Spectral properties of
the operator $H_0$}
    \label{ss21}
 Introduce the the annihilation
and the creation operators
$$
 a = a(b): = - 2i e^{-\varphi}
\frac{\partial}{\partial \overline{z}}\; e^{\varphi}, \; a^* =
a(b)^*: = - 2i e^{\varphi} \frac{\partial}{\partial z}\;
e^{-\varphi},
$$
 the function $\varphi$ being defined in \eqref{sn2}, and $z:
= x_1 + ix_2$, $\overline{z}:= x_1 - ix_2$. The operators $a$ and
$a^*$ defined initially on  $C_0^{\infty}(\rd)$, and then closed
in $L^2(\rd)$, are mutually adjoint. Set
    $$
    H_{\perp}^- = H_{\perp}^-(b) : = a^* a,  \quad H_{\perp}^+ = H_{\perp}^+(b) : = a a^*,
    $$
    $$
 H_{\perp} =
H_{\perp}(b) : = \left(
\begin{array} {cc}
H_{\perp}^- & 0\\
0 & H_{\perp}^+
\end{array}
\right) = H_{\perp}^- \oplus H_{\perp}^+.
    $$
 Then we have
    \bel{ker1}
    {\rm Ker}\;H_{\perp}^- = {\rm Ker}\;a = \left\{u \in
L^2(\rd) | u = ge^{-\varphi}, \;
    \frac{\partial g}{\partial
\overline{z}} = 0\right\},
    \ee
    $$
     {\rm Ker}\;H_{\perp}^+ =
{\rm Ker}\;a^* = \left\{u \in L^2(\rd) | u = ge^{\varphi}, \;
\frac{\partial g}{\partial z} = 0\right\},
    $$
    \bel{ker3}
     {\rm
Ker}\;H_{\perp} = \left\{{\bf u} = (u_1,u_2) | u_1 \in {\rm
Ker}\;H_{\perp}^-, \; u_2 \in {\rm Ker}\;H_{\perp}^+ \right\}.
    \ee
    Note that ${\rm
Ker}\;H_{\perp}^-$ (respectively, ${\rm Ker}\;H_{\perp}^+$) is a
weighted holomorphic (respectively, antiholomorphic) space of
Fock-Segal-Bargmann type (see e.g. \cite[Section 2 and Subsection
3.2]{H}). Since we have chosen $b_0
> 0$, and $\tilde{\varphi}$ is supposed to be bounded, we find
that ${\rm dim}\,{\rm Ker}\,H_{\perp}^- = \infty$ while ${\rm
dim}\,{\rm Ker}\,H_{\perp}^+ = 0$.
    \begin{pr} \label{p22} {\rm \cite[Proposition 1.2]{r5}}
Let $b$ be an  admissible magnetic field with $b_0 > 0$. Then $0
=\inf \sigma(H_{\perp})$ is an isolated eigenvalue of infinite
multiplicity. More precisely, we have
 \bel{sn32}
{\rm dim\,Ker}\;H_\perp = \infty,
    \ee
 and
   $$
     (0,{\mathcal C})
\subset \re \setminus \sigma(H_{\perp})
    $$
    with
    \bel{22}
    {\mathcal C}: =
2b_0 \exp{\left(-2 \; {\rm osc}\;\phe \right)},
    \ee
    where
${\rm osc}\;\phe: = \sup_{x \in \rd} \phe(x) - \inf_{x \in \rd}
\phe(x)$.
\end{pr}
{\em Remarks}: (i) Relation \eqref{sn32} holds true also for more
general backgrounds $\tilde{b}$. For example, it is sufficient
that $\tilde{b}$ is bounded, and the solution $ \tilde{\varphi}
\in C^2(\rd)$ of the Poisson equation \eqref{rt00} satisfies only
    \bel{sn37}
    \tilde{\varphi}(x) = o(|x|^2), \quad |x| \to \infty.
    \ee
    If $\tilde{b}$ is of the form \eqref{sn0}, and relations
    \eqref{sn33} - \eqref{sn35} (but not necessarily \eqref{sn36}) hold
    true, then
    $$
     \tilde{\varphi}(x)  : =  \int_{\rd} \frac{(\lambda \cdot x)^2}{|\lambda|^{2}}
     \int_0^1 (1-s) e^{is\lambda \cdot x} ds \,d\nu(\lambda), \quad x \in \rd,
    $$
    is in $C^2(\rd)$, and satisfies \eqref{rt00}
    and \eqref{sn37}.
    However, some of our further results, in particular,  Lemma \ref{gnl1}
    below, could be not true  for such more general magnetic
    fields.\\
 (ii) If $b$ is a periodic magnetic
field, the fact that the origin is an isolated eigenvalue of
$H_{\perp}$, was already mentioned in \cite{dno}, and was proved
in \cite{b}. A far going extension of the results of \cite{dno},
concerning the existence of {\em a strictly positive} isolated
eigenvalue of
$H_{\perp}$ of infinite multiplicity, could be found in \cite{vno}.\\

Now note that we have
    \bel{rt21}
    H_0^{\pm} = H_{\perp}^{\pm} \otimes I_{\parallel} + I_{\perp}
    \otimes H_{\parallel}
    \ee
    where $I_{\parallel}$ and $I_{\perp}$ are the identity
    operators in $L^2(\re)$ and $L^2(\rd)$ respectively, and
    $$
H_{\parallel} : = - \frac{d^2}{dx_3^2}
$$
is the self-adjoint operator, originally defined on
$C_0^{\infty}(\re)$, and then closed in $L^2(\re)$. Since
$\sigma(H_{\parallel})$ coincides  with $[0,\infty)$, and is
purely absolutely continuous, while $\inf \sigma(H_{\perp}^{-}) =
0$, we find that  \eqref{rt21} combined with, say, the arguments of \cite[Subsection 8.2.3]{abmbg}, implies the following
\begin{follow} \label{rtf21a}
Assume that $b$ is an admissible magnetic field. Then the spectrum
$\sigma(H_0)$ of the operator $H_0$ coincides with $[0,\infty)$,
and is purely absolutely continuous.
\end{follow}

\subsection{Berezin-Toeplitz operators}
    \label{ss22}
  Denote by $p = p(b)$ the orthogonal projection onto ${\rm
Ker}\,H_{\perp}^-(b)$ (see \eqref{ker1}). It is well known that
$p$ admits a continuous integral kernel ${\mathcal P}_b(x,y)$,
$x,y \in \rd$ (see e.g. \cite[Theorem 2.3]{H}).
    \begin{lemma} \label{gnl1}
    Assume that the magnetic field $b$ is admissible. Then we have
     \bel{rt50}
    \frac{b_0}{2\pi} e^{-2{\rm osc}\,\tilde{\varphi}} \leq
    {\mathcal P}_b(x,x) \leq \frac{b_0}{2\pi} e^{2{\rm
    osc}\,\tilde{\varphi}}, \quad x \in \rd.
    \ee
     \end{lemma}
    \begin{proof}
    Introduce the
functions \bel{rw} \phi_k(x): = \sqrt{\frac{b_0}{2\pi k!}}
\left(\frac{b_0}{2}\right)^{k/2} \left(x_1 + i x_2\right)^k
e^{-\varphi_0(x)}, \quad k \in {\mathbb Z}_+, \quad x \in \rd, \ee
which constitute an orthonormal in $L^2(\rd)$ basis of ${\rm
Ker}\;H_\perp^-(b_0) = {\rm Ker}\;a(b_0)$ (see e.g. \cite{rw}).
 Let  $\gamma
: l^2({\mathbb Z}_+) \to l^2({\mathbb Z}_+)$ be the operator given
in the canonic basis by the matrix
$\left\{g_{jk}\right\}_{j,k=0}^{\infty}$ with $g_{jk}: =
\int_{\rd} e^{-2 \phe} \phi_j \overline{\phi_k} \; dx$, $j,k \in
{\mathbb Z}_+$. It is easy to see that $\gamma$ is self-adjoint,
bounded, and
    \bel{above}
    \inf_{y \in \rd} e^{-2 \phe(y)} \leq \inf \sigma(\gamma)
    \leq \sup \sigma(\gamma) \leq \sup_{y \in \rd} e^{-2 \phe(y)}.
    \ee
    Set $\rho: = \gamma^{-1/2}$.
 Let $\left\{r_{jk}\right\}_{j,k=0}^{\infty}$ be the
matrix of  $\rho$ in the canonic basis of $l^2({\mathbb Z}_+)$.
Put
$$
\psi_j(x): = e^{-\,\phe(x)}\sum_{k=0}^{\infty} r_{jk} \phi_k(x),
\quad x \in \rd, \quad j \in {\mathbb Z}_+.
$$
Then $\left\{\psi_j\right\}_{j=0}^{\infty}$ is an orthonormal in
$L^2(\rd)$ basis of ${\rm Ker}\;a(b)$, and
    \bel{rrr}
    {\mathcal P}_{b}(x,x) =
    \sum_{j=0}^{\infty} |\psi_j(x)|^2 =
    e^{-2\,\phe(x)}\|\rho\phi(x)\|^2_{l^2({\mathbb Z}_+)}
    \ee
    where
$\phi(x): = \left\{\phi_k(x)\right\}_{k=0}^{\infty} \in
l^2({\mathbb Z}_+)$,  $x \in \rd$ being fixed (see \cite[Theorem
2.4]{H}). Making use of \eqref{above} and the spectral theorem, we
find that \eqref{rrr} and the obvious equality
$\sum_{k=0}^{\infty} |\phi_k(x)|^2 = \frac{b_0}{2\pi}$, valid for
each $x \in \rd$, imply
 \eqref{rt50}.
\end{proof}
The Berezin-Toeplitz operators necessary for the formulation of
our main results, have the form $p(b) U p(b)$ where $U : \rd \to
\re$. In Lemma \ref{frl0} below we describe a class of compact
operators of this type (admitting also complex-valued $U$).\\
 Let $X$ be a separable Hilbert space. In coherence with our
previous notations $S_1(X)$ and $S_2(X)$, we denote by $S_q(X)$,
$q \in [1,\infty)$, the Schatten - von Neumann classes of compact
linear operators $T$ for which the norm $\|T\|_q : = \left({\rm
Tr}\,|T|^q\right)^{1/q}$ is finite.
\begin{lemma} \label{frl0}
Let $U  \in L^q(\rd)$, $q \in [1,\infty)$. Assume that $b$ is an
admissible magnetic field. Then $p(b) U p(b) \in S_q(L^2(\rd))$,
and
    \bel{rt26}
    \|p(b) U p(b)\|_q^q \leq \frac{b_0}{2\pi} e^{2 {\rm
    osc}\,\tilde{\varphi}} \|U\|_{L^q}^q.
    \ee
\end{lemma}
\begin{proof}
If $U \in L^{\infty}(\rd)$, then
    \bel{rt51}
 \|p(b) U p(b)\| \leq \|U\|_{L^{\infty}}.
    \ee
    If $U \in L^1(\rd)$, then by $p(b) U p(b) = p(b) |U|^{1/2} e{^{i {\rm
    arg}\, U}} |U|^{1/2} p(b)$  and \eqref{rt50}, we have
    $$
    \|e{^{i {\rm
    arg}\, U}} |U|^{1/2} p(b)\|_2^2 =  \|p(b) |U|^{1/2}\|_2^2 = \int_{\rd}
    {\cal P}_b(x,x) |U(x)| dx \leq \frac{b_0}{2\pi} e^{2 {\rm
    osc}\,\tilde{\varphi}} \|U\|_{L^1}
    $$
    Therefore,
    \bel{rt52}
    \|p(b) U p(b)\|_1 \leq \frac{b_0}{2\pi} e^{2 {\rm
    osc}\,\tilde{\varphi}} \|U\|_{L^1}.
    \ee
Interpolating between \eqref{rt51} and \eqref{rt52}, we get
\eqref{rt26}.
\end{proof}
For further references, introduce the orthogonal projections
$$
P = P(b) : = p \otimes I_{\parallel}, \quad Q = Q(b): = I - P,
$$
acting in $L^2(\re^3)$, and the orthogonal projections
    \bel{fr53}
{\bf P} = {\bf P}(b) : = \left(
\begin{array} {cc}
P & 0\\
0 & 0
\end{array}
\right), \quad {\bf Q} = {\bf Q} (b) : = {\bf I} - {\bf P} =
 \left(
\begin{array} {cc}
Q & 0\\
0 & I
\end{array}
\right),
    \ee
acting in $L^2(\re^3; {\mathbb C}^2)$. Here $I$ and ${\bf I}$ are
the identity operators  in $L^2(\re^3)$ and  $L^2(\re^3; {\mathbb
C}^2)$ respectively.

\section{Main Results} \label{s3}
\setcounter{equation}{0} \subsection{Statement of the main
results} \label{ss23} For ${\bf x} = (x_1,x_2,x_3) \in \re^3$ we
denote by $\xp = (x_1,x_2)$ the variables on the plane
perpendicular to the magnetic field. Suppose that the matrix  $V$
satisfies
\begin{equation} \label{10a}
 v_{jk} \in C(\re^3), \quad  |v_{jk}({\bf
x})| \leq C_0  \langle \xp \rangle^{-\mper} \langle  x_3
\rangle^{-m_3}, \quad {\bf x} = (\xp, x_3) \in \re^3, \quad j,k =
1,2,
 \end{equation}
with $C_0 > 0$, $\mper >2$,  $m_3 >1$, and $\langle  y \rangle : =
(1 + |y|^2)^{1/2}$, $ y \in \re^d$, $d \geq 1$. Our main results
will be formulated under a more restrictive assumption than
\eqref{10a}, namely
\begin{equation} \label{fr1}
 v_{jk} \in C(\re^3), \quad |v_{jk}({\bf
x})| \leq C_0  \langle {\bf x} \rangle^{-m},  \quad {\bf x} \in
\re^3, \quad j,k = 1,2,
 \end{equation}
with $m > 3$. Note that \eqref{fr1} implies \eqref{10a} with any
$m_3 \in (0,m)$
and $\mper = m-m_3$.  \\
In the sequel we will assume that the perturbation of the operator
$H_0$ is of definite sign. For notational convenience, we will
suppose that
    \bel{rt22}
    V({\bf x}) \geq 0, \quad {\bf x} \in \rt,
    \ee
    and will consider the operators $H_0 + V$ or $H_0 - V$. \\
Assume that \eqref{10a} with $\mper > 2$, $m_3>1$, and
\eqref{rt22} hold true. Set
\begin{equation} \label{fr50}
W(\xp) : = \int_{\re} v_{11}(\xp,x_3) dx_3, \quad \xp \in \rd.
\end{equation}
If, moreover, $V$ satisfies \eqref{fr1}, then
\begin{equation} \label{fr19}
0 \leq W(\xp) \leq C_0'\langle \xp \rangle^{-m+1}, \quad \xp \in
\rd,
\end{equation}
where $C_0' = C_0 \int_{\re} \langle x \rangle^{-m} dx$. For $E >
0$ introduce the operator
\begin{equation} \label{fr15}
\omega(E) : = \frac{1}{2\sqrt{E}} p(b) W p(b).
 \end{equation}
Evidently, $\omega(E)$ is self-adjoint and non-negative in
$L^2(\rd)$. If $b$ is an admissible magnetic field, $E > 0$, and
$V$ satisfies \eqref{10a} with $m_{\perp} > 2$ and $m_3 > 1$, then
Lemma \ref{frl0} with $U = W$ implies $\omega(E) \in
S_1$.\\
Let $T= T^*$. Denote by ${\mathbb P}_{\delta}(T)$ the spectral
projection of $T$ associated with the Borel set $\delta \subset
\re$. Suppose that $T$ is compact and put
$$
n_{\pm}(s;T) : = {\rm rank}\,{\mathbb P}_{(s,\infty)}(\pm T),
\quad s>0.
$$
 Our first theorem concerns the asymptotic
behavior of the SSF $\xi(E;H,H_0)$ as the energy approaches the
origin from below.
\begin{theorem} \label{frt1}
Let \eqref{fr1} with $m>3$, and \eqref{rt22} hold true. Assume
that $b$ is an admissible magnetic field. Then for each
$\varepsilon \in (0,1)$ we have
\begin{equation} \label{fr17}
- n_+((1  - \varepsilon);\omega(E)) + O(1) \leq \xi( - E; H_0 - V,
H_0) \leq - n_+((1 + \varepsilon);\omega(E)) + O(1), \quad E
\downarrow 0.
 \end{equation}
\end{theorem}
{\em Remark}: By \eqref{101a}, if \eqref{rt22} holds true, then
$\xi(-E; H_0 + V, H_0) = 0$ for each $E > 0$. \\

Suppose again that the potential $V$ satisfies \eqref{10a} with
$\mper > 2$, $m_3>1$, and \eqref{rt22}. For $E
> 0$ define the matrix-valued function
\begin{equation} \label{fr51}
 {\cal W}_{E}(\xp) : = \left(
\begin{array} {cc}
w_{11}(\xp) & w_{12}(\xp)\\
w_{21}(\xp) & w_{22}(\xp)
\end{array}
\right), \quad  \xp \in \rd,
\end{equation}
where
$$
w_{11}(\xp) : = \int_{\re} v_{11}(\xp, x_3) \cos^2{(\sqrt{E}x_3)}
dx_3, \quad w_{22}(\xp) : = \int_{\re} v_{11}(\xp, x_3)
\sin^2{(\sqrt{E}x_3)} dx_3,
$$
$$
w_{12}(\xp) = w_{21}(\xp) : = \int_{\re} v_{11}(\xp, x_3)
\cos{(\sqrt{E}x_3)} \sin{(\sqrt{E}x_3)}dx_3.
$$
Set
\begin{equation} \label{fr17op}
\Omega(E) : = \frac{1}{2\sqrt{E}} p(b) {\cal W}_{E} p(b).
\end{equation}
Evidently, $\Omega(E)$ is self-adjoint  in $L^2(\rd; {\mathbb
C}^2)$, and $\Omega(E) \geq 0$. Since $\omega(E) \in S_1$, it is
easy to check that
$\Omega(E) \in S_1$ as well.\\
Our second theorem concerns the asymptotic behavior of the SSF
$\xi(E;H,H_0)$ as the energy  approaches the origin from above.
\begin{theorem} \label{frt2}
Let \eqref{fr1} with $m>3$ and \eqref{rt22} hold true. Assume that
$b$ is an admissible magnetic field. Then for each $\varepsilon
\in (0,1)$ we have
$$
\pm \frac{1}{\pi} {\rm Tr}\;{\rm arctan}\;((1  \pm
\varepsilon)^{-1}\Omega(E)) + O(1) \leq
$$
$$
\xi( E; H_0 \pm V, H_0) \leq
$$
\begin{equation} \label{fr18}
\pm \frac{1}{\pi} {\rm Tr}\;{\rm arctan}\;((1 \mp
\varepsilon)^{-1}\Omega(E)) + O(1), \quad E \downarrow 0.
 \end{equation}
\end{theorem}
    {\em Remark}: The privileged role of the entry $v_{11}$ of the
    matrix $V$ which occurs in the operators $\omega(E)$ and
    $\Omega(E)$, is determined by our assumption that  $b_0>0$, and,
    hence, the kernel of $H_{\perp}$ consists of elements with
    spin-up polarization (see \eqref{ker3}). In particular, we
    have
    $$
    {\bf P}(b) V {\bf P}(b) = \left(
    \begin{array} {cc}
    P(b) v_{11} P(b) & 0\\
    0 & 0
    \end{array}
    \right).
    $$
 The proofs of Theorems \ref{frt1} and \ref{frt2} can be found in
Section 4. In the following subsection we will describe explicitly
the asymptotics of $\xi( - E; H_0 - V, H_0)$ and $\xi( E; H_0 \pm
V H_0)$ as $E \downarrow 0$, under generic assumptions about the
behavior of $W(\xp)$ as $|\xp| \to \infty$.

 \subsection{Corollaries} \label{ss24}
  By \eqref{fr17} and \eqref{fr18},
we can reduce the analysis of the behavior as $E \to 0$ of $\xi(
E; H_0 \pm V,  H_0)$, to the investigation of the eigenvalue
asymptotics  of  compact Berezin-Toeplitz operators $p(b) U p(b)$,
discussed in the following three lemmas.\\
The first one treats the case where the decay of $U$ at infinity is
power-like. It involves the concept of {\em an integrated density of
states} (IDS) for the operator $H_{\perp}^-(b)$. Let $\chi_{Q}$ be
the characteristic function of the square $Q \subset \rd$, and let
$|Q|$ denote its area. We recall that the non-increasing function
$\varrho_b : \re \to [0,\infty)$ is called IDS for the operator
$H_{\perp}^-(b)$, if  it satisfies
    \bel{sn51}
\varrho_b(E) = \lim_{|Q| \to \infty} |Q|^{-1} {\rm Tr}\,(\chi_{Q}
{\mathbb P}_{(-\infty,E)} (H_{\perp}^-(b)) \chi_{Q})
    \ee
at its continuity points $E \in \re$ (see e.g. \cite{hlmw,
diwmi}). If $b=b_0$, i.e. if $\tilde{b}=0$, we have
    \bel{sn60}
    \varrho_{b_0}(E) = \frac{b_0}{2\pi} \sum_{q=0}^{\infty}
    \Theta(E-2b_0q), \quad E \in \re,
    \ee
    where $\Theta(t) = \left\{
    \begin{array} {l}
    0 \quad {\rm if} \quad t<0, \\
    1 \quad {\rm if} \quad t>0,
    \end{array} \right. $ is the Heaviside function.

\begin{lemma} \label{frl1} {\rm \cite[Proposition 3.5]{r5}}
Let  $U \in C^1(\rd)$ satisfy
$$
0 \leq U(\xp) \leq C_1 \langle \xp \rangle^{-\alpha},
\quad |\nabla U(\xp)| \leq C_1 \langle \xp \rangle^{-\alpha-1},
\quad \xp \in \rd,
$$
 for $\alpha > 0$ and $C_1 > 0$. Assume,
moreover, that:
\begin{itemize}
\item  $U(\xp) = u_0(\xp/|\xp|) |\xp|^{-\alpha}(1 + o(1))$ as
$|\xp| \to \infty$, where $u_0$ is a continuous function on
${\mathbb S}^1$ which does not vanish identically;
    \item $b$ is an admissible
magnetic field;
    \item there exists an  IDS $\varrho_b$ for the operator
$H_{\perp}^-(b)$.
\end{itemize}
Then we have
$$
n_+(s; p(b) U p(b)) = \frac{b_0}{2\pi} \left|\left\{\xp \in \rd |
U(\xp)
> s\right\}\right|\;(1 + o(1)) =
$$
    \bel{sn10}
 \Psi_{\alpha}(s; u_0,b_0)\; (1 +
o(1)), \quad s \downarrow 0,
    \ee
where, as above, $|.|$ denotes the Lebesgue measure, and
\begin{equation} \label{15d01}
\Psi_{\alpha}(s) = \Psi_{\alpha}(s; u_0,b_0): = s^{-2/\alpha}
\frac{b_0}{4\pi} \int_{{\mathbb S}^1}
u_0(\theta)^{2/\alpha}d\theta, \quad s>0.
\end{equation}
\end{lemma}
{\em Remarks}: (i) In \cite[Proposition 3.5]{r5} we considered only
the example of almost periodic  admissible magnetic fields, and
proved explicitly the existence of the IDS for the operator
$H_{\perp}^-(b)$. In Lemma \ref{frl1} above the existence of the IDS
is just a hypothesis. That is why,  we summarize here the main
ingredients of proof of \cite[Proposition 3.5]{r5} which do not
concern the existence of the IDS:
\begin{itemize}
    \item Applying variational and commutator techniques
    developed, in particular, in \cite{bpr, iwta}, we show that
    for each $E \in (0,{\mathcal C})$ we have
    \bel{sn11}
    n_+(s; p(b) U p(b)) = n_-(s/E; U^{1/2} (H_{\perp}^- - E)^{-1}
    U^{1/2}) (1 + o(1)), \quad s \downarrow 0;
    \ee
    \item Using the Birman-Schwinger principle, as well as the
    methods of \cite{adh, lev, hl} concerning the
    strong-electric-field asymptotics of the discrete spectrum of
    the operator $H_{\perp}^-(b) + g U$ lying in the gap $(0, {\mathcal
    C})$ of $\sigma(H_{\perp}^-(b))$, we obtain
    $$
    \lim_{g \to \infty} g^{-2/\alpha} n_-(g^{-1}; U^{1/2} (H_{\perp}^- - E)^{-1}
    U^{1/2}) =
    $$
    $$
    \int_{-\infty}^{E} \left|\left\{\xp  \in \rd \, |
    \,u_0(\xp/|\xp|) |\xp|^{-\alpha} > E-t\right\}\right|
    d\varrho_b(t) =
    $$
    \bel{sn12}
     E^{-2/\alpha} \frac{{\mathcal J}(b)}{2} \int_{{\mathbb S}^1} u_0(\theta)^{2/\alpha} d\theta, \quad E \in (0, {\mathcal
     C}),
     \ee
     where ${\mathcal J}(b)$ is the jump of the IDS $\varrho_b$ at
     the origin;
    \item We check that the family $H_\perp^-(b_0 + s\tilde{b})$, $s \in [0,1]$,
    is continuous in the norm resolvent sense, and, utilizing a gap-labelling theorem due to J. Bellissard
    \cite[Proposition 4.2.5]{bel2}, we find that the jump $J(b_0 +
    s\tilde{b})$ is independent of $s \in [0,1]$. In particular,
    \eqref{sn60} implies
    \bel{sn13}
    {\mathcal J}(b) = {\mathcal J}(b_0) = \frac{b_0}{2\pi}.
    \ee
    \end{itemize}
Putting together \eqref{sn11} -- \eqref{sn13}, we obtain
\eqref{sn10}. As a by-product of \eqref{sn51} with any $E \in
(0,{\mathcal C})$, and \eqref{sn13}, we obtain the formula
$$
\lim_{|Q| \to \infty} |Q|^{-1} \int_{Q} {\mathcal P}_b(\xp, \xp)
d\xp = \frac{b_0}{2\pi},
$$
valid if $b$ is an admissible magnetic field, and there exists an
IDS $\varrho_b$ for the operator $H_\perp^-(b)$.
\\
(ii) In the case $b = b_0$ (i.e. $\tilde{b} = 0$) a variant of Lemma
\ref{frl1} was proved in \cite{r0} with the help of
pseudo-differential techniques. In the case of general admissible
backgrounds $\tilde{b}$, the methods of \cite{r0} are not directly
applicable: due to the factor $\exp{(-\tilde{\varphi})}$ whose
derivatives generically do not decay at infinity, we do not obtain
suitable symbols of
pseudo-differential operators. \\

Our following two lemmas concern respectively the cases where $U$
decays exponentially at infinity, or has a compact support. First
note that, by \cite[Proposition 3.2]{r5}, we have
    \bel{sn5}
    n_+(\exp{(2\,{\rm osc}\tilde{\varphi})}s; p(b_0) U p(b_0)) \leq n_+(s; p(b) U p(b))
    \leq n_+(\exp{(-2\,{\rm osc}\tilde{\varphi})}s; p(b_0) U p(b_0)),
    \ee
    provided that $s>0$, $U : \rd \to [0,\infty)$, and the
    operator $U(-\Delta + 1)^{-1}$ is compact in $L^2(\rd)$.\\
 Combining \eqref{sn5} with the results of \cite[Proposition 3.1
with $q=0$]{rw} and of \cite[Proposition 3.2]{rw}, we obtain the
following
\begin{lemma} \label{frl2}
Let $0 \leq U \in L^{\infty}(\rd)$. Assume that
$$
\ln{U(\xp)} = -\eta |\xp|^{2\beta}(1 + o(1)), \quad |\xp| \to
\infty,
$$
for some $\beta \in (0,\infty)$,  $\eta \in (0,\infty)$. Let $b$
be an admissible magnetic field. Then we have
$$
n_+(s; p(b) U p(b)) = \Phi_{\beta}(s) (1 + o(1)), \quad s
\downarrow 0,
$$
where
\begin{equation} \label{15d02}
\Phi_{\beta}(s) = \Phi_{\beta}(s; \eta, b_0) : = \left\{
\begin{array} {l}
\frac{b_0}{2\eta^{1/\beta}} |\ln{s}|^{1/\beta} \quad {\rm if}
\quad 0 <
\beta < 1,\\
\frac{1}{\ln{(1+2\eta/b_0)}}|\ln{s}| \quad {\rm if} \quad
\beta = 1, \\
\frac{\beta}{\beta - 1}(\ln|\ln{s}|)^{-1}|\ln{s}| \quad {\rm if} \quad
1 < \beta < \infty,
\end{array}
\right. \quad s \in (0,e^{-1}).
\end{equation}
\end{lemma}
Similarly, the combination of \eqref{sn5} with  the results of
\cite[Proposition 3.2 with $q=0$]{rw} and of \cite[Proposition
3.2]{rw}, implies the following
\begin{lemma} \label{frl3}  Let $0 \leq U \in L^{\infty}(\rd)$. Assume that the
support of $U$ is compact, and that there exists a constant $C>0$
such that $U \geq C$ on an open non-empty subset of $\rd$. Let $b$
be an admissible magnetic field. Then we have
$$
n_+(s; p(b) U p(b)) =  \Phi_{\infty}(s) \; (1 + o(1)), \quad s
\downarrow 0,
$$
where
\begin{equation} \label{15d03}
\Phi_{\infty}(s) : = (\ln|\ln{s}|)^{-1}|\ln{s}|, \quad s \in
(0,e^{-1}).
\end{equation}
\end{lemma}
Employing now Lemma \ref{frl1}, \ref{frl2},  or \ref{frl3},  we
find that \eqref{fr17} immediately entails the following
\begin{follow} \label{frf01}
Let \eqref{fr1} with $m > 3$, and \eqref{rt22}  hold true.\\
(i)  Assume that the hypotheses of Lemma \ref{frl1} hold with $U =
W$ and $\alpha = m-1$. Then
$$
\xi(- E; H_0-V,H_0) = - \frac{b_0}{2\pi} \left|\left\{\xp \in \rd
| W(\xp) > 2\sqrt{E}\right\}\right|\;(1 + o(1)) =
$$
\begin{equation} \label{fr80}
- \Psi_{m-1}(2\sqrt{E}; u_0,b_0) \; (1 + o(1)), \quad E
  \downarrow 0,
\end{equation}
the function $\Psi_{\alpha}$ being defined in \eqref{15d01}. \\
(ii) Assume that the hypotheses of Lemma \ref{frl2} hold with $U =
W$. Then we have
$$
\xi(- E; H_0-V,H_0) = - \Phi_{\beta}(2\sqrt{E};\eta,b_0) \; (1 +
o(1)), \quad E \downarrow 0, \quad \beta \in (0,\infty),
$$
the functions $\Phi_{\beta}$ being defined in \eqref{15d02}.\\
(iii) Assume that the hypotheses of Lemma \ref{frl3} hold with $U
= W$. Then we have
$$
\xi( - E; H_0-V,H_0) = - \Phi_{\infty}(2{\sqrt{E}}) \; (1 + o(1)),
\quad E \downarrow 0,
$$
the function $\Phi_{\infty}$ being defined in \eqref{15d03}.
\end{follow}
    {\em Remark}: By \eqref{101a}, the results of Corollary
    \ref{frf01}, as well those of Theorem \ref{frt1}, concern the
    asymptotic distribution near the origin of the (negative)
    discrete spectrum of the operator $H_0 - V$. Results, related
    to Corollary \ref{frf01} (i) concerning perturbations $V$ of power-like decay,
    could be found in \cite{iwta}
    where, similarly to the present article, magnetic fields ${\bf
    B} = (0,0,b)$ of constant direction are considered.
    Moreover, in \cite{iwta}, the perturbation $V$ is not obliged to be
    asymptotically homogeneous, the decay rate $m$ is allowed to
    be any positive number, and two distinct types of asymptotic
    formulae concerning the case $m \in (0,2)$ and $m \in
    (2,\infty)$ are deduced, the latter being similar to
    \eqref{fr80}. On the other hand, in \cite{iwta} the function $b$ is assumed to be positive,
    its derivative is supposed to decay at infinity, and the perturbation $V$ is scalar.
    Results which extend Lemma \ref{frl3}, and are related to
    Corollary \ref{frf01} (iii), are contained in \cite{fp}.\\

Next, the combination of Theorem \ref{frt2} with Lemmas \ref{frl1}
-- \ref{frl3} yields the following
\begin{follow} \label{frf02}
(i) Let \eqref{fr1} with $m > 3$, and \eqref{rt22}, hold true.
Assume that the hypotheses of Lemma \ref{frl1} are fulfilled for
$U = W$ and $\alpha = m-1$. Then we have
$$
\xi( E; H_0 \pm V,H_0) = \pm \frac{b_0}{2\pi^2} \int_{\rd} {\rm
arctan}\;((2\sqrt{E})^{-1}W(\xp)) d\xp \; (1 + o(1)) =
$$
$$
\pm  \, \frac{1}{2 \cos{(\pi/(m-1))}}\,
\Psi_{m-1}(2\sqrt{E};u_0,b_0) \,(1+o(1)) , \quad E \downarrow 0.
$$
(ii) Let \eqref{fr1} with $m > 3$, and \eqref{rt22}, hold true.
Suppose in addition that $V$ satisfies \eqref{10a}  for some
$m_{\perp} > 2$ and $m_3
> 2$. Finally, assume that the hypotheses of Lemma \ref{frl2} are
fulfilled for $U = W$. Then we have
$$
\xi(E; H_0 \pm V,H_0) = \pm \, \frac{1}{2} \,
\Phi_{\beta}(2\sqrt{E};\eta,b_0) \; (1 + o(1)), \quad E \downarrow
0, \quad \beta \in (0,\infty).
$$
(iii) Let the assumptions of the previous part be fulfilled,
except that the hypotheses of  Lemma \ref{frl2} are replaced by
those of Lemma \ref{frl3}. Then we have
$$
\xi(E; H_0 \pm V, H_0) = \pm \, \frac{1}{2} \,
\Phi_{\infty}(2\sqrt{E}) \; (1 + o(1)), \quad E \downarrow 0.
$$
\end{follow}
The main ingredient of the proof of Corollary \ref{frf02} is the
estimate
    \bel{sn20}
    {\rm Tr}\,\arctan{(s^{-1} \Omega(E))} = {\rm Tr}\,\arctan{(s^{-1}
    \tilde{\Omega}(E))}(1 + o(1)), \quad E \downarrow 0, \quad s>0,
    \ee
    where
    $$
 \tilde{\Omega}(E): = \frac{1}{2\sqrt{E}} p(b) \left(
 \begin{array} {cc}
 W & 0\\
 0 & 0
 \end{array}
 \right) p(b), \quad E > 0,
 $$
 $W$ being defined in \eqref{fr50}.
 Estimate \eqref{sn20} is obtained by using the Lifshits-Krein trace
 formula \eqref{sn15} with $f(E) = \arctan{E}$, $E \in \re$.
Since the argument of the proof  of Corollary \ref{frf02} is
completely analogous to the one of \cite[Corollary 3.2]{fr}, we
omit the details.\\

{\em Remark}: By \eqref{sn30}, Corollary \ref{frf02} as well as
Theorem \ref{frt2} concern the low-energy asymptotics of the
scattering phase ${\rm arg}\, {\rm det}\,S(H_0 \pm V,H_0)$.\\

Putting together the results of Corollaries \ref{frf01} and
\ref{frf02} for negative perturbations, we obtain
\begin{follow} \label{frf03}
Under the assumptions of Corollary \ref{frf02} (i) we have
    \bel{rt60}
    \lim_{E \downarrow 0} \frac{\xi(E; H_0 - V, H_0)}{\xi(-E; H_0 - V, H_0)} = \frac{1}{2\cos(\pi/(m-1))}, \quad m>3,
    \ee
    while under the assumptions of Corollary \ref{frf02} (ii)--(iii) we have
    \bel{rt61}
    \lim_{E \downarrow 0} \frac{\xi(E; H_0 - V, H_0)}{\xi(-E; H_0 - V, H_0)} = \frac{1}{2}.
    \ee
    \end{follow}
    {\em Remark}: Formulae \eqref{rt60} -- \eqref{rt61} could be interpreted
    as {\em generalized Levinson formulae}.
    We recall that the classical Levinson formula relates the (finite) limiting values as $E \uparrow 0$ and
    $E \downarrow 0$ of the SSF $\xi(E; -\Delta + V; -\Delta)$ where $\Delta$ is the Laplacian in $\re^d$,
    $d \geq 1$, and $V : \re^d \to \re$ is a scalar potential which decays fast enough at infinity
    (see the original work \cite{l} or the survey article \cite{rob}).

\section{Auxiliary results}
\setcounter{equation}{0} \label{s4}

\subsection{A representation of the SSF}
\label{ss31}
 In
this subsection we introduce a suitable representation of the SSF
$\xi(E; H_0 \pm V, H_0)$, $E \in (-\infty,{\mathcal
C})\setminus\{0\}$, based on a general  abstract result of A.
Pushnitski  \cite{p1}. \\
Assume that $V$ satisfies \eqref{rt22} and \eqref{10a}. Set
    \bel{rt23}
L({\bf x}) = \left\{\ell_{jk}({\bf x})\right\}_{j,k=1}^2 : =
V({\bf x})^{1/2}, \quad {\bf x} \in \re^3.
    \ee
Then for $E < 0$ we have
    \bel{36}
    L(H_0 - E)^{-1/2} \in S_{\infty}(L^2(\rd; {\mathbb C}^2)),
    \ee
    \bel{37}
    L(H_0 - E)^{-1} \in S_{2}(L^2(\rd; {\mathbb C}^2)).
    \ee
 For $z \in
{\mathbb C}_+ : = \left\{\zeta \in {\mathbb C}\, | \, {\rm
Im}\,\zeta > 0\right\}$, set $T(z): = L (H_{0} - z)^{-1} L$. By
\cite{be} (see also \cite[Lemma 4.1]{p1}), for almost every $E \in
\re$ the operator-norm limit
\begin{equation} \label{38}
T(E + i0): = {\rm n}-\lim_{\delta \downarrow 0} T(E + i\delta)
\end{equation}
exists, and
\begin{equation} \label{39}
{\rm Im}\,T(E + i0) \in S_1.
\end{equation}
For trivial reasons the limit in \eqref{38} exists, and \eqref{39}
holds for {\em each} $E < 0 = \inf \sigma(H_0)$. In Corollary
\ref{fff} below we show that this is also true for each $E \in (0,
{\mathcal C})$. Hence, by \cite[Lemma 2.1]{p1}, the quantity
\begin{equation} \label{39dnes}
\tilde{\xi}(E ; H_0 \pm V, H_0) = \pm \int_{\re} n_{\mp}(1; {\rm
Re} \, T(E + i0) + t\; {\rm Im} \, T(E + i0))\; d\mu (t), \quad E
\in (-\infty,{\mathcal C}) \setminus \{0\},
\end{equation}
with
$$
d\mu(t) : = \frac{dt}{\pi(1+t^2)},
$$
 is well-defined. Arguing as in the proof of
\cite[Proposition 2.5]{bpr} (see also \cite[Proposition
2.1]{brs}), and bearing in mind Proposition \ref{p31}, Corollary
\ref{f31}, and Proposition \ref{frp1} below, we easily prove the
following
    \begin{pr} \label{rtf21}
    Assume that $V$ satisfies \eqref{10a} with $m_\perp > 2$, $m_3 > 1$, and \eqref{rt22}.
    Suppose that $b$ is an admissible magnetic filed. Then
    $\tilde{\xi}(\cdot; H_0 \pm V, H_0)$ is bounded on every
    compact subset of $(-\infty, {\mathcal C})\setminus \{0\}$, and
    is continuous on $(-\infty, {\mathcal C})\setminus (\{0\} \cup \sigma_{\rm pp}(H \pm V))$.
    \end{pr}
Since $V$ satisfies \eqref{10a} with $m_\perp > 2$, $m_3 > 1$,
relation \eqref{13} holds true and the SSF $\xi(E; H_0 \pm V,
H_0)$ is well defined for almost every $E \in \re$. On the other
hand, by \cite[Theorem 1.2]{p1} we have
    $$
\xi(E ; H_0 \pm V, H_0) = \tilde{\xi}(E ; H_0 \pm V, H_0)
    $$
for almost every $E \in \re$. In this article we identify $\xi(E ;
H_0 \pm V, H_0)$ with $\tilde{\xi}(E ; H_0 \pm V, H_0)$
for $E \in (-\infty,{\mathcal C}) \setminus \{0\}$.\\

{\em Remark}: The representation of the SSF described  above
 admits a generalization to non-sign-definite perturbations
$V$ (see \cite{gm, p3}). This generalization is based on the
concept of the index of orthogonal projections (see \cite{ass}).\\
We formulate our main results and their corollaries for the case
of
    perturbations of constant sign because certain key auxiliary
    facts are known to be true only in this case.

\subsection{Estimates of sandwiched resolvents} \label{ss32} For $z \in {\mathbb
C}_+ $
 define the operator
 $R(z): = \left( - \frac{d^2}{dx_3^2} - z\right)^{-1}$, bounded in
 $L^2(\re)$.
The operator $R(z)$ admits the integral kernel ${\cal R}_z(x_3 -
x_3')$ where
  ${\cal R}_z(x) = i e^{i\sqrt{z} |x|}/(2\sqrt{z})$,   $x \in \re$,
  and the branch of $\sqrt{z}$ is  chosen so that ${\rm Im}\;\sqrt{z}
  > 0$.\\
 For $z \in {\mathbb
C}_+ $
 introduce the operators
    \bel{rt55}
T_<(z) : = L {\bf P}(H_0 - z)^{-1} L,  \quad T_>(z) : = L {\bf
Q}(H_0 - z)^{-1} L,
    \ee
bounded in  $L^2(\rt; {\mathbb C}^2)$ (see \eqref{fr53} for the
definition of the
  orthogonal projections ${\bf P}$ and ${\bf Q}$).  Then we have  $T_<(z) =
  L
\;\Big( (p \otimes R(z)) \oplus 0) \Big) \; L$. \\
For $E \in \re$, $E \neq  0$, define  $R(E)$
  as
the operator with integral kernel ${\cal R}_{E}(x_3-x_3')$ where
\begin{equation} \label{vg1}
{\cal R}_E(x) : = \lim_{\delta \downarrow 0}{\cal R}_{E +
i\delta}(x) = \left\{
\begin{array} {l}
\frac{e^{-\sqrt{-E} |x|}}{2\sqrt{-E}} \quad {\rm if} \quad
E < 0,\\
\frac{i e^{i\sqrt{E} |x|}}{2\sqrt{E}} \quad {\rm if} \quad E
> 0,
\end{array}
\right.
\quad x \in \re.
\end{equation}
For $E \in \re$, $E \neq 0$,  set
$$
T_<(E) : = L \;\Big(  (p \otimes R(E)) \oplus 0\Big) \; L.
$$
\begin{pr} \label{p31}
Let \eqref{10a} with $m_\perp>2$, $m_3>1$, and \eqref{rt22} hold
true. Then the operator-valued function $\overline{{\mathbb
C}_+}\setminus \{0\} \ni z \mapsto T_<(z) \in S_1$ is well defined
and continuous. Moreover,
\begin{equation} \label{311}
\|T_<(E)\|_1 \leq C_1 (1 + E_+^{1/4})|E|^{-1/2}, \quad E \in \re
\setminus \{0\},
\end{equation}
with $C_1$ independent of $E$.
\end{pr}
\begin{proof}
  The operator
 $T_<(z)$ admits the representation
\begin{equation} \label{pu9}
T_<(z) = M \;((G \otimes J(z))\oplus 0) \; M, \quad z \in
\overline{{\mathbb C}_+}\setminus \{0\},
\end{equation}
where $M: L^2(\rt; {\mathbb C}^2) \to L^2(\rt; {\mathbb C}^2)$ is
the multiplier by the matrix-valued function
    \bel{rt74}
M(\xp, x_3) : = \langle \xp
\rangle^{\mper/2} \langle x_3 \rangle^{m_3/2} L(\xp,x_3), \quad
(\xp,x_3) \in \rt,
    \ee
the operator $G : = \langle \xp \rangle^{-\mper/2} p \,\langle \xp
\rangle^{-\mper/2}$ acts in $L^2(\rd)$, while  $$J(z) : = \langle
x_3 \rangle^{-m_3/2} R(z) \langle x_3 \rangle^{-m_3/2}$$ acts in
$L^2(\re)$. Evidently,
    \bel{rt56}
\|T_<(z)\|_1 \leq \|M\|^2 \|G\|_1 \|J(z)\|_1, \quad z \in
\overline{{\mathbb C}_+}\setminus \{0\}.
    \ee
By \eqref{10a}, the operator $M$ is bounded. Further, $\|G\|_1 =
\|pUp\|_1$ with $U(\xp) = \langle \xp \rangle^{-\mper}$, $\xp \in
\rd$. By $\mper > 2$ we have $U \in L^1(\rd)$, and  Lemma
\ref{frl0} implies $G \in S_1$. Moreover, $M$ and $G$ are
independent of $z$. By \cite[Subsection 4.1]{bpr} the
operator-valued function $\overline{{\mathbb C}_+}\setminus \{0\}
\ni z \to J(z) \in S_1$ is well defined and continuous, and admits
the estimate
\begin{equation} \label{rt57}
\|J(E)\|_1 \leq C_1' (1 + E_+^{1/4})|E|^{-1/2}, \quad E \in \re
\setminus \{0\},
\end{equation}
with $C_1'$ independent of $E$. Now the claim of the lemma follows
from \eqref{pu9} -- \eqref{rt57}.
\end{proof}
For further references we state here the following obvious
 \begin{follow} \label{f31}
Let $V$ satisfy the assumptions of Proposition \ref{p31}. Let $E
\in \re$,   $E \neq 0$. Then ${\rm Im}\;T_<(E) \geq 0$. Moreover,
if $E < 0$, then ${\rm Im}\;T_<(E) = 0$.
\end{follow}
\begin{pr} \label{frp1}
Let $V$ satisfy the assumptions of Proposition \ref{p31}. Then the
function ${\mathbb C}\setminus [{\mathcal C}, \infty) \ni z \mapsto
T_>(z) \in S_2$ is well defined and analytic. Moreover, for $E \in
(-\infty, {\mathcal C})$ we have
    \bel{rev1}
T_>(E) = T_>(E)^*,
    \ee
    and
    \bel{rt73}
    \|T_>(E)\|_2 \leq C_2\left(1 + \frac{(E+1)_+}{{\mathcal C} - E}\right),
    \ee
    with $C_2$ independent of $E$.
\end{pr}
\begin{proof}
We have
$$
T_>(z) = L((Q(H_0^- -z)^{-1}) \oplus (H_0^+ -z)^{-1})L, \quad z \in {\mathbb C}\setminus[{\mathcal C},\infty).
$$
The function ${\mathbb C}\setminus [{\mathcal C}, \infty) \ni z
\mapsto T_>(z) \in {\mathcal B}$,  the class of linear bounded
operators, is well defined and analytic, and \eqref{rev1} holds true
for $E \in (-\infty, {\mathcal C})$,  just because ${\mathbb
C}\setminus [{\mathcal C}, \infty)$ is included in the resolvent
sets of the operator $H_0^-$ defined on $Q D(H_0^-)$, and of the
operator $H_0^+$ defined on $D(H_0^+)$. Further, set
$$
F(\xp, x_3) = \langle \xp \rangle^{-\mper/2} \langle x_3 \rangle^{-m_3/2}, \quad (\xp, x_3) \in \rt.
$$
Note that $L = FM$, the matrix $M$ being defined in \eqref{rt74}. Then we have
    \bel{rt70}
    \|T_>(z)\|_2^2 \leq  \|L\|^2 \left( \|Q (H_0^- -z)^{-1}F\|_2^2 + \|(H_0^+ -z)^{-1}F\|_2^2\right)
    \|M\|^2.
    \ee
     Applying the spectral theorem for bounded functions of self-adjoint operators, the resolvent identity,
     and the diamagnetic inequality for Hilbert-Schmidt operators, we get
    $$
    \|Q (H_0^- -z)^{-1}F\|_2 \leq C(z) \|(H_0^- + 1)^{-1}F\|_2 \leq
    $$
    $$
    C(z) \|1 + (H_0^- + 1)^{-1} b\| \|((i\nabla + {\bf A})^2 + 1)^{-1}F\|_2 \leq
    $$
    \bel{rt71}
    C(z) \|1 + (H_0^- + 1)^{-1} b\| \|(-\Delta + 1)^{-1}F\|_2
    \ee
    where
    $$
    C(z) : = \sup_{s \in [{\mathcal C},\infty)} \left|\frac{s + 1}{s - z}\right|,
    \quad z \in {\mathbb C}\setminus [{\mathcal C}, \infty).
    $$
    Similarly,
    \bel{rt72}
      \|(H_0^+ -z)^{-1}F\|_2 \leq C(z) \|1 - (H_0^+ + 1)^{-1} b\| \|(-\Delta + 1)^{-1}F\|_2.
    \ee
Since $\|(-\Delta + 1)^{-1}F\|_2 < \infty$, we find that \eqref{rt70} -- \eqref{rt72} imply that
$T_>(z) \in S_2$ if $z \in {\mathbb C}\setminus[{\mathcal C}, \infty)$, and that \eqref{rt73} holds true. \\
The analyticity of $T_>(z)$ in $S_2$ follows from an appropriate estimate of the Hilbert-Schmidt norm
of the derivative $\frac{dT_>(z)}{dz}$.
 \end{proof}
 Propositions \ref{p31} and \ref{frp1} immediately entail
\begin{follow} \label{fff}
Let $V$ satisfy the assumptions of Proposition \ref{p31}. Then for
$E = (-\infty,{\mathcal C})\setminus\{0\}$ the operator-norm limit
\eqref{38} exists, and
\begin{equation} \label{fff1}
T(E+i0) = T_<(E) + T_>(E).
\end{equation}
Moreover,
\begin{equation} \label{fff3}
{\rm Re}\;T(E+i0) = {\rm Re}\;T_<(E)  + T_>(E),
\end{equation}
\begin{equation} \label{fff2}
{\rm Im}\;T(E+i0) = {\rm Im}\;T_<(E).
\end{equation}
\end{follow}

 \section{Proof of  the main results }
\setcounter{equation}{0}
    \label{s5}
\subsection{A preliminary estimate}
    \label{ss41}
This subsection contains a preliminary estimate (see \eqref{fr3}
below) which will be used in the proofs of  Theorems \ref{frt1} --
\ref{frt2}. \\

The following lemma contains a suitable version of the Weyl
inequalities for the eigenvalues of compact operators.

\begin{lemma} \label{snl1} {\rm \cite[Chapter I, Eq. (1.32)]{bs}}
Let $T_j^*$, $j=1,2$, be compact self-adjoint operators acting in
the same Hilbert space. Then we have
    \bel{wineq}
    n_{\pm}(s_1 + s_2; T_1 + T_2) \leq n_{\pm}(s_1; T_1) + n_{\pm}(s_2;
    T_2)
    \ee
    for every $s_1> 0$ and $s_2 > 0$.
    \end{lemma}

\begin{pr} \label{frp1a}
Let \eqref{10a} with $m>3$, and  \eqref{rt22} hold true. Let $E =
(-\infty,{\mathcal C})\setminus\{0\}$. Then the asymptotic
estimates
$$
\int_{\re} n_{\pm}(1+\varepsilon; {\rm Re} \, T_<(E) + t\; {\rm
Im} \, T_<(E))\; d\mu (t) + O(1) \leq
$$
$$
\int_{\re} n_{\pm}(1; {\rm Re} \, T(E + i0) + t\; {\rm Im} \, T(E
+ i0))\; d\mu (t)  \leq
$$
\begin{equation} \label{fr3}
\int_{\re} n_{\pm}(1-\varepsilon; {\rm Re}\, T_<(E) + t\; {\rm Im}
\, T_<(E))\; d\mu (t) + O(1)
\end{equation}
hold as $E \to 0$ for each $\varepsilon \in (0,1)$.
\end{pr}
\begin{proof}
By \eqref{fff3} and \eqref{fff2}, and the Weyl inequalities
\eqref{wineq}, we have
$$
\int_{\re} n_{\pm}(1+\varepsilon; {\rm Re} \, T_<(E) + t\; {\rm
Im} \, T_<(E))\; d\mu (t) -  n_{\mp}(\varepsilon; T_>(E)) \leq
$$
$$
\int_{\re} n_{\pm}(1; {\rm Re} \, T(E + i0) + t\; {\rm Im} \, T(E + i0))\; d\mu (t)  \leq
$$
\begin{equation} \label{fr3a}
\int_{\re} n_{\pm}(1-\varepsilon; {\rm Re}\, T_<(E) + t\; {\rm Im}
\, T_<(E))\; d\mu (t) +  n_{\pm}(\varepsilon;  T_>(E)).
\end{equation}
 Evidently,
 $n_{\pm}(\varepsilon;  T_>(E)) \leq \varepsilon^{-2}
 \|T_>(E)\|_2^2$,
which combined with \eqref{rt73}, yields
    \bel{rt58}
 n_{\pm}(\varepsilon;  T_>(E)) = O(1), \quad E \to 0.
 \ee
 Now \eqref{fr3} follows from \eqref{fr3a} and \eqref{rt58}.
 \end{proof}

\subsection{Proof of Theorem \ref{frt1}}
\label{ss42}
 Throughout the subsection we
assume  the hypotheses of Theorem \ref{frt1}.
 By Corollary \ref{f31} we have ${\rm
Im}\;T_{<}( - E) = 0$ and, hence, ${\rm Re}\;T_<( - E) = T_<(-
 E)$ if $E > 0$. Therefore,
\begin{equation} \label{fr4}
\int_{\re} n_{+}(s; {\rm Re}\, T_<(-E) + t\; {\rm Im} \, T_<(-E))\;
d\mu (t) = n_{+}(s;  T_<(-E)), \quad E>0, \quad s>0.
\end{equation}
For  $E > 0$ define ${\cal O}(E) : L^2(\re^3; {\mathbb C}^2) \to
L^2(\re^3; {\mathbb C}^2)$ as the operator with matrix-valued
integral kernel
$$
\frac{1}{2\sqrt{E}} \ell_{j1}(\xp,x_3)\;{\cal
P}_b(\xp,\xp')\ell_{1k}(\xp',x_3'), \quad j,k =1,2, \quad
(\xp,x_3), (\xp',x_3') \in \re^3.
$$
\begin{pr} \label{frp2}
For each $\varepsilon \in (0,1)$ and $s>0$ we have
\begin{equation} \label{fr6}
n_{+}((1+\varepsilon)s ; {\cal O}(E)) + O(1) \leq n_{+}(s;  T_<( -
E)) \leq n_{+}((1-\varepsilon)s ; {\cal O}(E)) + O(1), \quad E
\downarrow 0.
\end{equation}
\end{pr}
\begin{proof}
Fix $s>0$ and $\varepsilon \in (0,1)$. By the Weyl inequalities
\eqref{wineq},
$$
n_{+}((1+\varepsilon)s ; {\cal O}(E)) - n_-(\varepsilon s ; T_<( -
E) - {\cal O}(E)) \leq
$$
$$
 n_{+}(s;  T_<( -
E)) \leq
$$
$$
n_{+}((1-\varepsilon)s ; {\cal O}(E)) + n_+(\varepsilon s ; T_<( -
E) - {\cal O}(E)).
$$
In order to get \eqref{fr6}, it suffices to show that there exists
a compact operator $\tilde{T}$ such that
\begin{equation} \label{fr23}
{\rm n}-\lim_{E \downarrow 0} (T_<( - E) - {\cal O}(E)) =
\tilde{T}.
\end{equation}
Pick $m' \in (3,m)$, and note that
$$
T_<( - E) - {\cal O}(E) = \tilde{M}_{m,m'} ((\tilde{G}_{m-m'}
\otimes \tilde{J}_{m'}(E)) \oplus 0) \tilde{M}_{m,m'}
$$
where $\tilde{M}_{m,m'}$ is the multiplier by the bounded
matrix-valued function
$$
\langle \xp \rangle^{(m-m')/2} \langle
x_3 \rangle^{m'/2}L(\xp,x_3), \quad (\xp,x_3) \in \re^3,$$
$\tilde{G}_{m-m'} : L^2(\rd) \to  L^2(\rd)$ is the operator with
integral kernel
$$
\langle \xp \rangle^{-(m-m')/2} {\cal P}_b(\xp,\xp')\langle \xp'
\rangle^{-(m-m')/2}, \quad \xp, \xp' \in \rd,
$$
and $\tilde{J}_{m'}(E)$, $E > 0$, is the operator with integral
kernel
$$
-\frac{1}{2\sqrt{E}}\langle x_3
\rangle^{-m'/2}\left(1-e^{-\sqrt{E}|x_3-x_3'|}\right) \langle x_3'
\rangle^{-m'/2}, \quad x_3, x_3' \in \re.
$$
 Set
 \bel{511a}
\tilde{T}= \tilde{M}_{m,m'} ((\tilde{G}_{m-m'} \otimes
\tilde{J}_{m'}(0)) \oplus 0) \tilde{M}_{m,m'}
    \ee
where $\tilde{J}_{m'}(0) : L^2(\re) \to L^2(\re)$ is the operator
with integral kernel
    $$
    -\frac{1}{2}\langle x_3
\rangle^{-m'/2} |x_3-x_3'| \langle x_3'  \rangle^{-m'/2}, \quad
x_3, x_3' \in \re.
$$
Note that $\tilde{T}$ admits a matrix-valued integral kernel
    \bel{511b}
-\frac{1}{2}\ell_{j1}(\xp,x_3) |x_3-x_3'|{\cal
P}_b(\xp,\xp')\ell_{1k}(\xp',x_3'), \; j,k=1,2, \;(\xp,x_3),
\;(\xp',x_3') \in \re^3.
    \ee
 Since $m-m' > 0$, the operator $\tilde{G}_{m-m'}$ is compact
by Lemma \ref{frl0}. Since $m' > 3$ we have $\tilde{J}_{m'}(E) \in
S_2$ for $E\geq0$. Bearing in mind that $\tilde{M}_{m,m'}$ is
bounded, we find that the operator $\tilde{T}$ is compact.
Finally, we have $\lim_{E \downarrow 0}\|\tilde{J}_{m'}(E) -
\tilde{J}_{m'}(0)\|_2 = 0$ which easily implies \eqref{fr23}.
\end{proof}
\begin{pr} \label{frp3}
For each
 $E > 0$, and $s>0$ we have
\begin{equation} \label{fr7}
n_{+}(s ; {\cal O}(E)) = n_{+}(s ; \omega(E)),
\end{equation}
the operator $\omega(E)$ being defined in \eqref{fr15}.
\end{pr}
\begin{proof}
Define the operator $K : L^2(\rt; {\mathbb C}^2) \to  L^2(\rd)$ by
$$
(K {\bf u})(\xp) : = \sum_{k=1,2}\int_{\rd} \int_{\re} {\cal
P}_b(\xp, \xp') \ell_{1k}(\xp',x_3') u_k(\xp',x_3')\, dx_3'\,
d\xp', \quad \xp \in \rd,
$$
where ${\bf u} = \left( \begin{array} {c} u_1\\u_2 \end{array}
\right) \in L^2(\rt; \cd)$. We have
$$
{\cal O}(E) = \frac{1}{2\sqrt{E}} K^* K, \quad {\omega}(E) =
\frac{1}{2\sqrt{E}} K\,K^*.
$$
Since $n_+(s; K^* K) = n_+(s; K\,K^*)$ for each $s>0$, we get \eqref{fr7}.
\end{proof}
Putting together  \eqref{39dnes}, \eqref{fr3},
 \eqref{fr4}, \eqref{fr6},
and \eqref{fr7}, we get  \eqref{fr17}, which concludes the proof
of Theorem \ref{frt1}.

\subsection{Proof of Theorem \ref{frt2}}
\label{ss43}
 Throughout the subsection we
assume  the hypotheses of Theorem \ref{frt2}.
\begin{pr} \label{frp4}
 For each $s>0$ we have
\begin{equation} \label{fr8}
n_{\pm}(s;  {\rm Re}\,T_<(E)) = O(1), \quad E \downarrow 0.
\end{equation}
\end{pr}
\begin{proof}
The operator ${\rm Re}\,T_<(E)$ admits the matrix-valued integral
kernel
$$
-\ell_{j1}(\xp,x_3) \frac{\sin{(\sqrt{E}|x_3-x_3'|)}}{2\sqrt{E}}
{\cal P}_b(\xp,\xp')\ell_{1k}(\xp',x_3'), \; j,k=1,2, \;
(\xp,x_3), \,(\xp',x_3') \in \re^3.
$$
 Arguing as
in the proof of Proposition \ref{frp2}, we find that ${\rm
n}-\lim_{E \downarrow 0}  {\rm Re}\,T_<( E) = \tilde{T}$ (see
\eqref{511a} --  \eqref{511b}) which implies \eqref{fr8}.
\end{proof}
Making use of Propositions \ref{frp1a} and \ref{frp4} and Corollary
\ref{f31}, as well as of the Weyl inequalities \eqref{wineq} and the
evident identities
$$
\int_{\re} n_{\pm}(s; t T) d\mu(t) = \frac{1}{\pi} {\rm Tr}\;{\rm
arctan}\;(s^{-1} T), \quad s>0,
$$
with $T = T^* \geq 0$, $T \in S_1$, we obtain the following
\begin{follow} \label{frf2}
For each $\varepsilon \in (0,1)$ and $s>0$ we have
$$
\frac{1}{\pi} {\rm Tr}\;{\rm arctan}\;((s(1+\varepsilon))^{-1}
{\rm Im}\;T_<(E)) + O(1) \leq
$$
$$
\int_{\re} n_{\pm}(s; {\rm Re}\;T_<(E) + t \;{\rm Im}\;T_<(E))
d\mu(t) \leq
$$
\begin{equation} \label{fr9}
\frac{1}{\pi} {\rm Tr}\;{\rm arctan}\;((s(1-\varepsilon))^{-1}
{\rm Im}\;T_<(E)) + O(1), \quad E \downarrow 0.
\end{equation}
\end{follow}
\begin{pr} \label{frp5}
For each
 $E > 0$ and $s>0$ we have
\begin{equation} \label{fr20}
n_{+}(s ; {\rm Im}\;T_<(E)) = n_{+}(s ; \Omega(E)),
\end{equation}
the operator $\Omega(E)$ being defined in \eqref{fr17op}.
Consequently,
\begin{equation} \label{fr21}
{\rm Tr}\;{\rm arctan}\;(s^{-1}{\rm Im}\;T_<(E)) = {\rm Tr}\;{\rm
arctan}\;(s^{-1}\Omega(E)), \quad E>0, \quad s > 0.
\end{equation}
\end{pr}
\begin{proof}
The operator ${\rm Im}\,T_<(E)$ admits the matrix-valued integral
kernel
$$
\ell_{j1}(\xp,x_3) \frac{\cos{(\sqrt{E}(x_3-x_3'))}}{2\sqrt{E}}
{\cal P}_b(\xp,\xp')\ell_{1k}(\xp',x_3'), \; j,k=1,2, \;
(\xp,x_3), \, (\xp',x_3') \in \re^3.
$$
Define the operator ${\cal K} : L^2(\rt; \cd) \to  L^2(\rd; \cd)$
by
$$
{\cal K}{\bf u} : = {\bf v} = \left( \begin{array} {c} v_1 \\ v_2
\end{array} \right) \in L^2(\rd; \cd), \quad {\bf u} = \left( \begin{array} {c} u_1 \\ u_2
\end{array} \right)\in L^2(\rt; \cd),
$$
where
$$
v_1(\xp): = \sum_{k=1,2}\int_{\rd} \int_{\re} {\cal P}_b(\xp,
\xp') \cos({\sqrt{E}x_3'}) \ell_{1k}(\xp',x_3') u_k(\xp',x_3')\,
dx_3'\, d\xp',
$$
$$
v_2(\xp): = \sum_{k=1,2}\int_{\rd} \int_{\re} {\cal P}_b(\xp,
\xp') \sin({\sqrt{E}x_3'}) \ell_{1k}(\xp',x_3') u_k(\xp',x_3')\,
dx_3'\, d\xp', \quad \xp \in \rd.
$$
We have
$$
{\rm Im}\;T_<(E) = \frac{1}{2\sqrt{E}} {\cal K}^* {\cal K}, \quad
\Omega(E)  = \frac{1}{2\sqrt{E}} {\cal K}\,{\cal K}^*.
$$
Since $n_+(s; {\cal K}^* {\cal K}) = n_+(s; {\cal K}\,{\cal K}^*)$ for each $s>0$, we get \eqref{fr20}.
\end{proof}
Now the combination of \eqref{39dnes}, \eqref{fr3}, \eqref{fr9},
and \eqref{fr21} yields \eqref{fr18}.\\

{\large \bf Acknowledgements}. The author thanks Rafael Tiedra de
Aldecoa for an illuminating discussion on the contents  of
\cite{abmbg}.\\ The partial support of the Chilean Science
Foundation {\em Fondecyt} under Grant 1090467, and of {\em N\'ucleo
Cient\'ifico ICM} P07-027-F ``{\em Mathematical Theory of Quantum
and Classical Magnetic Systems"}, is gratefully acknowledged.

\end{document}